March 03, 2015

# Poisson's formula with principal value integrals and some special Gradshteyn and Ryzhik integrals

**Khristo N. Boyadzhiev**

Department of Mathematics and Statistics,

Ohio Northern University, Ada, OH 45810, USA

k-boyadzhiev@onu.edu

**Abstract.** Poisson's integral formula for holomorphic functions on the right half plane can be used to quickly evaluate certain integrals from Gradshteyn and Ryzhik's table. In addition we prove a version of Poisson's formula for principal value integrals and use it in several interesting cases.

## 1. Introduction

Let $f(z)$, $z = x + iy$, be a bounded analytic function on the open right half plane $x > 0$ (RHP). For such functions Poisson's integral formula holds (see [14], chapter 6, with some adjustment)

(1.1) $$f(x+iy) = \frac{1}{\pi} \int_{-\infty}^{+\infty} f(it) \frac{x}{x^2 + (y-t)^2} \, dt ,$$

where $x > 0$, $-\infty < y < \infty$ and $f(it)$ stands for the boundary values of the function (they exist almost everywhere). For real valued integrals we shall use the two equations

(1.2) $$\operatorname{Re} f(x+iy) = \frac{x}{\pi}\int_{-\infty}^{+\infty} \operatorname{Re} f(it)\frac{dt}{x^2+(y-t)^2}$$

(1.3) $$\operatorname{Im} f(x+iy) = \frac{x}{\pi}\int_{-\infty}^{+\infty} \operatorname{Im} f(it)\frac{dt}{x^2+(y-t)^2}$$

This formula has numerous applications in harmonic analysis, potential theory, and partial differential equations. Books on complex variables present examples of integral evaluation based usually on the residue theorem or on Cauchy's integral formula. We shall demonstrate that in certain cases Poisson's formula is a better tool for integral evaluation.

One such example was given in Section 8 in [2]. In the present paper we shall use Poisson's formula to evaluate a number of integrals from Gadshteyn and Ryzhik's table [6]. In Section 3 we prove a version of this formula for principal value integrals and apply it to several special cases. Section 5 contains some historical remarks on the works of Hardy on principal value integrals.

## 2. Examples

**Example 1**. Let $f(z) = e^{-az} = e^{-ax}e^{-iay} = e^{-ax}(\cos ay - i\sin ay)$, where $a > 0$ and $z = x + iy$. The function is analytic and bounded for $x > 0$. Therefore,

(2.1) $$e^{-ax}(\cos ay - i\sin ay) = \frac{x}{\pi}\int_{-\infty}^{+\infty}\frac{\cos at - i\sin at}{x^2+(y-t)^2}dt \ .$$

Separating real and imaginary parts we find

(2.2) $$\int_{-\infty}^{+\infty}\frac{\cos at}{x^2+(y-t)^2}dt = \frac{\pi}{x}e^{-ax}\cos ay \ , \quad \int_{-\infty}^{+\infty}\frac{\sin at}{x^2+(y-t)^2}dt = \frac{\pi}{x}e^{-ax}\sin ay$$

Usually the evaluation of these integrals is not so straightforward. Evaluation by residues, in fact, is equivalent to the proof of Poisson's formula.

**Example 2**. We shall apply the same method to the function $f(z) = ze^{-az}$ which is not bounded on vertical lines in the RHP. To avoid this obstacle we assume that Poisson's formula is applied to the bounded on the RHP function $f(z) = \dfrac{ze^{-az}}{1+\varepsilon z}$ for small $\varepsilon > 0$ and then, in the final result, we set $\varepsilon \to 0$. As long as the integrals are at least conditionally convergent (as they are), the formula works. Thus we get

$$(2.3) \qquad e^{-ax}(x\cos ay + y\sin ay) = \frac{x}{\pi} \int_{-\infty}^{+\infty} \frac{t\sin at}{x^2 + (y-t)^2} \, dt,$$

$$(2.4) \qquad e^{-ax}(y\cos ay - x\sin ay) = \frac{x}{\pi} \int_{-\infty}^{+\infty} \frac{t\cos at}{x^2 + (y-t)^2} \, dt.$$

In particular, with $y = 0$ in (2.3) we have

$$(2.5) \qquad e^{-ax} = \frac{1}{\pi} \int_{-\infty}^{+\infty} \frac{t\sin at}{x^2 + t^2} \, dt = \frac{2}{\pi} \int_{0}^{+\infty} \frac{t\sin at}{x^2 + t^2} \, dt$$

which is entry 3.723.3 in Gradshteyn and Ryzhik's table [6]. With $y = 0$ the first integral in (2.2) becomes

$$(2.6) \qquad \int_{-\infty}^{+\infty} \frac{\cos at}{x^2 + t^2} \, dt = \frac{\pi}{x} e^{-ax}$$

(entry 3.723.2 in that table). A combination of (2.5) and (2.6) with $x = 1$ gives

$$(2.7) \qquad \int_{0}^{+\infty} \frac{\cos at + t\sin at}{1 + t^2} \, dt = \pi e^{-a}$$

which is 3.784.6 in [6].

**Example 3.** Here we use the function $f(z) = e^{-bz} \cos az$, where $0 < a < b$. Clearly,

$$(2.8) \qquad f(it) = \cos bt \cos at - i \sin bt \cos at.$$

We apply Poisson's formula, separate the real part, and set $y = 0$ to get

$$(2.9) \qquad \int_{0}^{+\infty} \frac{\cos bt \cos at}{x^2 + t^2} \, dt = \frac{\pi}{2x} e^{-bx} \cosh ax,$$

as the integrand is even. This is entry 3.742.3. in [6] in a slightly different form. The result is true also for $a = b$.

**Example 4**. Now we shall prove entry 3.742.5 in [6] which includes two cases. Taking $f(z) = ze^{-bz} \cos az$, $0 < a < b$, we use the same $\varepsilon$-trick as in Example 2. We compute $\operatorname{Re} f(it) = t \cos at \sin bt$, and then with $y = 0$,

$$(2.10) \qquad \int_0^{+\infty} \frac{t \sin bt \cos at}{x^2 + t^2} dt = \frac{\pi}{2} e^{-bx} \cosh ax$$

(as the integrand is even). Changing the function to $f(z) = ze^{-bz} \sinh az$, this time $\operatorname{Re} f(it) = -t \cos bt \sin at$, and we find

$$(2.11) \qquad \int_0^{+\infty} \frac{t \cos bt \sin at}{x^2 + t^2} dt = \frac{-\pi}{2} e^{-bx} \sinh ax .$$

Adding (2.10) to (2.11) yields

$$(2.12) \qquad \int_0^{+\infty} \frac{t \sin(a+b)t}{x^2 + t^2} dt = \frac{\pi}{2} e^{-(a+b)x},$$

which is, in fact, (2.5). Note that (2.10) and (2.11) can be obtained directly from (2.9) by differentiation for $b$ and $a$ correspondingly.

**Example 5**. Let $\alpha > 0$ and consider the function $f(z) = \dfrac{\log(1+\alpha z)}{z}$ with the principal value of the logarithm. Then $f(it) = \dfrac{-i}{t} \log(1 + i\alpha t)$ and

$$(2.13) \qquad \operatorname{Re} f(it) = \frac{1}{t} \operatorname{Arg}(1 + i\alpha t) = \frac{\arctan \alpha t}{t} .$$

With $y = 0$ we find from (1.1)

$$(2.14) \qquad \frac{\log(1+\alpha x)}{x^2} = \frac{1}{\pi} \int_{-\infty}^{+\infty} \frac{\arctan \alpha t}{t} \frac{dt}{x^2 + t^2} = \frac{2}{\pi} \int_0^{+\infty} \frac{\arctan \alpha t}{t} \frac{dt}{x^2 + t^2}$$

which is entry 4.535.9 in [6].

At the end of this section we present two more examples with integrals from [6]. The functions to which Poisson's formula is applied are aso listed. Details are left to the reader.

**Example 6.** This is entry 3.725 (3) from [6]. Prove that for every $x > 0$,

$$(2.15) \qquad \int_0^{+\infty} \frac{\sin at \cos bt}{t} \frac{dt}{x^2 + t^2} = \frac{\pi}{2x^2} e^{-bx} \sinh ax, \text{ when } 0 < a < b,$$

Here $f(z) = \frac{1}{z} e^{-bz} \sinh az$. Also in this entry,

$$(2.16) \qquad \int_0^{+\infty} \frac{\sin at \cos bt}{t} \frac{dt}{x^2 + t^2} = \frac{\pi}{2x^2} (1 - e^{-ax} \cosh bx), \text{ when } 0 < b < a.$$

For this case we use the function $f(z) = \frac{1}{z}(1 - e^{-az} \cosh bz)$.

**Example 7.** Prove that for every $x > 0, a > 0, b > 0$,

$$(2.17) \qquad \int_0^{+\infty} \frac{\cos bt}{b^2 + t^2} \frac{dt}{x^2 + t^2} = \frac{\pi}{2bx} \frac{be^{-ax} - xe^{-ab}}{b^2 - x^2}$$

(entry 3.728 in [6]). Hint: Use the function $f(z) = \dfrac{be^{-az} - ze^{-ab}}{b^2 - z^2}$.

### 3. Poisson's formula for principal value integrals

We can evaluate by Poisson's formula also principal value integrals.

**Example 8.** Taking $f(z) = \operatorname{sech} az, a > 0$, with $z = x + iy$, we find from (1.2)

$$(3.1) \qquad P\int_{-\infty}^{+\infty} \frac{1}{\cos at} \frac{dt}{x^2 + (y-t)^2} = \operatorname{Re} \frac{\pi}{x \cosh az},$$

since $\cosh(ait) = \cos at$. Setting $y = 0$ and using the fact that the integrand is even, we find

$$(3.2) \quad P\int_0^{+\infty} \frac{1}{\cos at} \frac{dt}{x^2+t^2} = \frac{\pi}{2x\cosh ax}$$

This evaluation is justified by the proposition below (a modification of Poisson's formula). As we shall see, the imaginary part of $f(z)$ is represented by a certain series (not appearing here).

Before stating the proposition we need some definitions. Recall that if the function $f(x)$ is defined on $[a,b]$ except at some point $c$, where $a<c<b$, then

$$(3.3) \quad P\int_a^b f(x)\,dx = \lim_{\varepsilon \to 0}\left(\int_a^{c-\varepsilon} f(x)\,dx + \int_{c+\varepsilon}^b f(x)\,dx\right).$$

The function $\sec at$ in example (3.1) has an infinite sequence of simple poles on the real axis. We shall give below a formal definition of a principal value integral for functions of this nature; and this definition will be used in Poisson's formula.

**Definition 1**. Let $\{ic_k, c_k \in \mathbb{R}\}$ be a finite or infinite sequence of numbers on the $y$-axis such that $|c_k - c_m| \geq d > 0, k \neq m$. For every $0 < \varepsilon < d/2$ we define the contour $L(\varepsilon)$ to be the $y$-axis indented to the right at the points of this sequence, so that each $ic_k$ is isolated from the RHP by a small semicircle centered at $ic_k$ and with radius $\varepsilon$. Here $\varepsilon$ is so small that the semicircles do not touch each other. The contour is oriented upward. Next, we define the set $M(\varepsilon)$ to be $L(\varepsilon)$ without the small semicircles, that is, $M(\varepsilon)$ is the union of all segments $[ic_{k-1} + i\varepsilon, ic_k - i\varepsilon]$ on the y-axis. When the sequence is finite, say, $ic_0,...,ic_m$, we add to $M(\varepsilon)$ also the infinite intervals $(-i\infty, ic_0 - i\varepsilon]$ and $[ic_m + i\varepsilon, +i\infty)$. Likewise, when the sequence is bounded from one side only, we add one such infinite interval at that side.

**Definition 2.** Let $ic_k$ and $M(\varepsilon)$ be as in the above definition. If the function $g(z)$ is defined on the imaginary axis except at the points $ic_k$, then we define

$$(3.4) \quad P\int_{-\infty}^{\infty} g(it)\,dt = \lim_{\varepsilon \to 0} \int_{M(\varepsilon)} g(it)\,dt\ .$$

**Proposition**. *Suppose $f(z)$ is a holomorphic function on the closed RHP with simple poles $ic_k$ on the $y$-axis, which are separated as in Definition 1. We want the residues at these poles to be purely imaginary, that is, of the form $ie_k$, with $e_k$ real numbers. Suppose $\varepsilon > 0$, and let $L(\varepsilon)$ be the indented contour defined for the sequence of poles. We assume that the function $f(z)$ is bounded to the right of $L(\varepsilon)$ (including on $L(\varepsilon)$ itself) for every $\varepsilon > 0$. Finally, we require that the series*

$$(3.5) \qquad \sum_k \frac{e_k}{x^2 + (c_k - y)^2}$$

*be convergent for any $x > 0$, $y$ real. In that case, Poison's formula (1.2) holds and the integral exists as a principal value integral. That is, for any $x > 0$ and any $-\infty < y < \infty$ we have*

$$(3.6) \qquad \operatorname{Re} f(x+iy) = \frac{1}{\pi} P \int_{-\infty}^{+\infty} \operatorname{Re} f(it) \frac{x}{x^2 + (y-t)^2} dt .$$

Proof. We shall mimic the classical proof of Poisson's formula [14, chapter 6] for bounded functions on the RHP. For any small $\varepsilon > 0$ we consider the contour $L(\varepsilon)$ defined above. Then we take a closed, counterclockwise oriented contour $C(R)$ consisting of a semicircle in the RHP centered at the origin, with radius $R$, and with ends $\pm iR$ lying on $L(\varepsilon)$ between the poles. On the left side $C(R)$ is closed by the part of $L(\varepsilon)$ between $\pm iR$. For a given $z = x+iy$ with $x > 0$ we take $R$ so large that $z$ is inside $C(R)$. By Cauchy's integral formula,

$$(3.7) \qquad f(z) = \frac{1}{2\pi i} \oint_{C(R)} f(\lambda) \frac{d\lambda}{\lambda - z} .$$

At the same time, for the point $z^* = -x+iy$ which is outside $C(R)$ we have

$$(3.8) \qquad 0 = \frac{1}{2\pi i} \oint_{C(R)} f(\lambda) \frac{d\lambda}{\lambda - z^*} .$$

By subtraction

(3.9) $$f(z) = \frac{1}{2\pi i} \oint_{C(R)} f(\lambda) \left( \frac{1}{\lambda - z} - \frac{1}{\lambda - z^*} \right) d\lambda = \frac{1}{\pi i} \oint_{C(R)} f(\lambda) \frac{x}{(\lambda - iy)^2 - x^2} d\lambda .$$

When $R \to \infty$ the integral on the proper semicircle approaches zero and the integration, after adjusting the direction, reduces to $L(\varepsilon)$

(3.10) $$f(z) = \frac{-1}{\pi i} \int_{L(\varepsilon)} f(\lambda) \frac{x}{(\lambda - iy)^2 - x^2} d\lambda = \frac{1}{\pi i} \int_{L(\varepsilon)} f(\lambda) \frac{x}{x^2 - (\lambda - iy)^2} d\lambda$$

Now we look at the small semicircles separating the poles. Let $C_k(\varepsilon)$ be one such semicircle with radius $\varepsilon$ and centered at the simple pole $ic_k$ where the residue is $ie_k$, and $c_k, e_k$, are real numbers. Then by Jordan's lemma (Theorem 1 in 3.1.4, [15])

(3.11) $$\lim_{\varepsilon \to 0} \frac{1}{\pi i} \int_{C_k(\varepsilon)} f(\lambda) \frac{x}{x^2 - (\lambda - iy)^2} d\lambda = \frac{ie_k x}{x^2 + (c_k - y)^2} .$$

Note that this value is purely imaginary. From (3.9)

(3.12) $$f(z) = \lim_{\varepsilon \to 0} \frac{1}{\pi i} \int_{M(\varepsilon)} f(\lambda) \frac{x}{x^2 - (\lambda - iy)^2} d\lambda + \sum_k \frac{ie_k x}{x^2 + (c_k - y)^2}$$

and with $\lambda = it$ in this integral we have

(3.13) $$f(z) = \lim_{\varepsilon \to 0} \frac{1}{\pi} \int_{M(\varepsilon)} f(it) \frac{x}{x^2 + (t - y)^2} dt + \sum_k \frac{ie_k x}{x^2 + (c_k - y)^2} .$$

The limit here is the principal value integral. Separating the real parts here we find

(3.14) $$\operatorname{Re} f(z) = \frac{1}{\pi} \int_{-\infty}^{\infty} \operatorname{Re} f(it) \frac{x}{x^2 + (t - y)^2} dt .$$

The proof is completed.

The function $f(z) = \operatorname{sech} az$ in Example 8 has simple poles $ic_k = (2k + 1)\frac{\pi i}{2a}$, $k = 0, \pm 1, \pm 2, \ldots$, with residues $ie_k = i(-1)^{k+1}$. The corresponding series (3.5) is obviously convergent.

## 4. Examples with principal value integrals

The collection of Gradshteyn and Ryzhik [6] contains several very interesting principal value integrals which can be evaluated by the modified Poisson formula.

**Example 9**. For every $a > 0$,

$$(4.1) \qquad P\int_0^{+\infty} \frac{t}{\sin at} \frac{dt}{x^2 + t^2} = \frac{\pi}{2\sinh ax}$$

This is Gradshteyn and Ryzhik's integral 3.747.3, similar to (3.1). To prove it we apply Poisson's formula to $f(z) = \frac{z}{\sinh az}$, $a > 0$. Since the function is unbounded on vertical lines, we use again the $\varepsilon$-trick from Example 2 (i.e. we work with $f(z) = \frac{z}{(1+\varepsilon z)\sinh az}$ first, and then set $\varepsilon \to 0$), This yields

$$(4.2) \qquad \operatorname{Re} \frac{z}{\sinh az} = \frac{1}{\pi} P\int_{-\infty}^{+\infty} \frac{t}{\sin at} \frac{x}{x^2 + (y-t)^2} dt$$

Setting $y = 0$ and then reducing by $x$, we come to (4.1), as the integrand is even.

Note that here the function $f(z)$ has simple poles $\frac{k\pi i}{a}$, $k = 0, \pm 1, \pm 2, \ldots$, with residues $\frac{(-1)^k k\pi i}{a^2}$.

**Example 10**. We shall present here a group of five integrals extending five Gradshteyn and Ryzhik entries, namely, 3.743 (1-4) and 3.744, the last one in the list below. They all can be proved by Poisson's formula (3.6). The function $f(z)$, $z = x + iy$, for which the formula is applied, is listed on the LHS of the equations. The entries in [6] correspond to $y = 0$. In each case $0 < a < b$.

$$(4.3) \qquad \operatorname{Re} \frac{\sinh az}{\sinh bz} = \frac{1}{\pi} P\int_{-\infty}^{+\infty} \frac{\sin at}{\sin bt} \frac{x}{x^2 + (y-t)^2} dt,$$

$$\text{(4.4)} \qquad \operatorname{Re}\frac{z\sinh az}{\cosh bz} = \frac{-1}{\pi} P\int_{-\infty}^{+\infty}\frac{t\sin at}{\cos bt}\frac{x}{x^2+(y-t)^2}\,dt,$$

$$\text{(4.5)} \qquad \operatorname{Re}\frac{z\cosh az}{\sinh bz} = \frac{1}{\pi} P\int_{-\infty}^{+\infty}\frac{t\cos at}{\sin bt}\frac{x}{x^2+(y-t)^2}\,dt,$$

$$\text{(4.6)} \qquad \operatorname{Re}\frac{\cosh az}{\cosh bz} = \frac{1}{\pi} P\int_{-\infty}^{+\infty}\frac{\cos at}{\cos bt}\frac{x}{x^2+(y-t)^2}\,dt,$$

$$\text{(4.7)} \qquad \operatorname{Re}\frac{\sinh az}{z\cosh bz} = \frac{1}{\pi} P\int_{-\infty}^{+\infty}\frac{\sin at}{t\cos bt}\frac{x}{x^2+(y-t)^2}\,dt.$$

Note that when $y=0$ all integrands here are even functions in $t$ and the integrals on the entire line $(-\infty,\infty)$ equal twice the integrals on $[0,\infty)$. All these integrals appear also in table [5] as Fourier sine or cosine transforms. When $y=0$, we can reduce by $x$ in equations (4.4) and (4.5), so the integrals become correspondingly

$$\text{(4.8)} \qquad \frac{\sinh ax}{\cosh bx} = \frac{-1}{\pi}\operatorname{PV}\int_{-\infty}^{+\infty}\frac{\sin at}{\cos bt}\frac{t}{x^2+t^2}\,dt = \frac{-2}{\pi}\operatorname{PV}\int_{0}^{+\infty}\frac{\sin at}{\cos bt}\frac{t}{x^2+t^2}\,dt,$$

$$\text{(4.9)} \qquad \frac{\cosh ax}{\sinh bx} = \frac{1}{\pi}\operatorname{PV}\int_{-\infty}^{+\infty}\frac{\cos at}{\sin bt}\frac{t}{x^2+t^2}\,dt = \frac{2}{\pi}\operatorname{PV}\int_{0}^{+\infty}\frac{\cos at}{\sin bt}\frac{t}{x^2+t^2}\,dt.$$

The above examples, we hope, illustrate well the method of integral evaluation by Poisson's formula. There are many other integrals in Gadshteyn and Ryzhik's table [6] (and also in other places) that can be approached by this technique.

**Remark.** In [3] several hyperbolic integrals similar to some integrals in the present paper were evaluated. For example, the following integral was evaluated

$$\text{(4.13)} \qquad \frac{1}{\pi}\int_{0}^{+\infty}\frac{1}{\cosh at}\frac{x}{x^2+t^2}\,dt = 2\sum_{k=1}^{\infty}\frac{(-1)^{k-1}}{2ax+(2k-1)\pi},$$

This is equation (7.1) in [3] and entry 3.522.3 in [6]. The integral is very similar to (3.2), with a hyperbolic cosine instead of trigonometric cosine. There are, however, big differences, as the integral in (4.13) is uniformly convergent, while the one in (3.2) is a principal value integral. The Poison formula cannot be applied for (4.13) in a straightforward manner. If we consider the function $f(z) = \dfrac{1}{\cos az}$ on the RHP, then its boundary value on the $y$-axis is indeed $f(it) = \dfrac{1}{\cosh at}$, but this function has a sequence of poles in the RHP and Poisson's formula is not applicable. We need a different method to compute explicitly the harmonic function

$$(4.14) \qquad u(x+iy) = \frac{1}{\pi} \int_{-\infty}^{+\infty} \frac{1}{\cosh at} \frac{x}{x^2 + (y-t)^2} \, dt \ .$$

## 5. Some historical remarks

Principal value integrals were introduced by Cauchy, who evaluated many such integals by the residue method. In particular, Cauchy evaluated in his book [4] the five integrals (4.3) – (4.7) (see the notes on p.121 in [7] and also paragraph 6 on p. 373 in [8]). His solutions were long and elaborate, but the new theory and the method were revolutionary. Many years after Cauchy, at the beginning of the 20$^{th}$ century Hardy undertook a thorough study of principal value integrals, providing a detailed and strict definition, proving various properties, and evaluating many integrals of the type presented here. Hardy was obviously very much involved in this matter, as he published six papers dedicated to principal value integrals in the period 1900 -1909, namely [7] – [12]. Hardy's papers on integral calculus can be found in the fifth volume of his collected works [13]. A contemporary treatment of such principal value integrals by the method of contour integration is given by Antimirov et al in chapter 8 of the book [1].

In [7], pp. 137-139, Hardy evaluated the integral

$$(5.1) \qquad P\int_0^{+\infty} \cos at \tan bt \, \frac{t \, dt}{x^2 + t^2} = \frac{\pi \cosh ax}{e^{2bx} + 1} \ ,$$

where $0 \le a < b$. This integral can be proved by Poisson's formula using the function

$$(5.2) \qquad f(z) = \frac{z \cosh az}{e^{2bz} + 1} = \frac{z e^{-bz} \cosh az}{2 \cosh bz} \ .$$

Hardy showed that the integral is continuous for $a$ at $a = 0$ and therefore,

(5.3) $$P\int_0^{+\infty} \tan bt \, \frac{t}{x^2+t^2} \, dt = \frac{\pi}{e^{2bx}+1}$$

which is entry 3.749 (1) in [6]. Following Hardy ([10], pp. 85-86) we can make the following observation about this integral: Putting $a = b$ in (4.8) we arrive at

(5.4) $$P\int_0^{+\infty} \tan bt \, \frac{t}{x^2+t^2} \, dt = \frac{-\pi}{2} \tanh bx$$

which contradicts (5.3). This shows that the integral in (4.8) is not continuous for $a$ at $a = b$. A similar remark can be made for the integral in (4.9).

**References**

[1]     M. Ya. Antimirov, Andrei A. Kolyshkin and Rémi Vaillancour, Complex Variables, Academic Press, 1998.

[2]     Khristo N, Boyadzhiev, Louis Medina, and Victor Moll, The integrals in Gradshteyn and Ryzhik. Part 11: The incomplete beta function, SCIENTIA: Series A: Mathematical Sciences, 18 (2009), 61-75. (see section 8)

[3]     Khristo N. Boyadzhiev, Victor Moll, The integrals in Gradshteyn and Ryzhik.Part 21: Hyperbolic functions, SCIENTIA: Series A: Mathematical Sciences, v. 22 (2012), 109-127.

[4]     Augustin-Louis Cauchy, Mémoire sur les intégrales définies, 1827, Œuvres complètes, série 1, tome 1, 319-506.

[5]     Arthur Erdélyi et al, Tables of Integral Transforms, Volume 1, McGrow-Hill, 1954.

[6]     Izrail S. Gradshteyn and Iosif M. Ryzhik, Table of Integrals, Series, and Products, Academic Press, 1980.

[7]     Godfrey H. Hardy, On differentiation and integration under the integral sign, Quarterly J. Pure and Applied Math, 32 (1901), 66-140. (Continued as Quarterly J. Math.)

[8]     Godfrey H. Hardy, General theorems in contour integration with some applications, Quarterly J. Pure and Applied Math, 32 (1901), 369-384.

[9]     Godfrey H. Hardy, The elementary theory of Cauchy's principal values. Proc. London Math. Soc. (1), 34, (1902), 16-40.

[10]    Godfrey H. Hardy, The theory of Cauchy's principal values. (Second Paper: the use of principal values in some of the double limit problems of the integral calculus.), Proc. London Math. Soc. (1), 34, (1902), 55-91.

[11]    Godfrey H. Hardy, The theory of Cauchy's principal values (Third Paper: differentiation and integration of principal values), Proc. London Math. Soc. (1), 35, (1903), 81-107.

[12]    Godfrey H. Hardy, The theory of Cauchy's principal values (Fourth paper : The integration of principal values – continued…), Proc. London Math. Soc. (2), 7, (1909), 181-208.

[13]    Godfrey H. Hardy, Collected papers of G. H. Hardy, volume V, Oxford, 1972.

[14]    Paul Koosis, Introduction to $H_p$ spaces, Cambridge University Press, 1980.

[15]    Dragoslav S. Mitrinović and Jovan D. Keckic, The Cauchy Method of Residues: Theory and Applications, Sringer 2001.